\documentclass[a4paper,11pt]{article}

\usepackage{amsmath}
\usepackage{amsfonts}
\usepackage{amssymb}
\usepackage{caption}
\usepackage{graphicx}
\usepackage{graphics}
\usepackage{color}
\usepackage{ushort}
\usepackage{makeidx}
\usepackage{float}
\usepackage{appendix}
\usepackage{accents}
\usepackage{framed}
\usepackage{versions}
\usepackage{xcolor}
\usepackage{mathbbol}
\usepackage{fontenc}
\usepackage[utf8]{inputenc}
\usepackage{blindtext}
\usepackage{emptypage}
\usepackage{verbatim}
\usepackage{enumitem}
\usepackage{array}
\usepackage{multirow}
\usepackage{extramarks}

\def\hmath$#1${\texorpdfstring{{\rmfamily\textit{#1}}}{#1}}

\usepackage[bookmarks]{hyperref}
\hypersetup{
    colorlinks=true,   	
    linkcolor=red,      
    citecolor = [rgb]{0 0.7 0},   	
    filecolor=magenta, 	
    urlcolor=blue
}

\newcommand{\equald}{\mbox{$ \;\stackrel{\cal D}{=}\; $}}

\newcommand{\RL}{{\mathbb R}}

\newcommand{\IN}{{\mathbb Z}}

\newcommand{\VAR}{\mbox{\rm Var}}

\def\ba{\begin{align}}
\def\ea{\end{align}}
\def\ban{\begin{align*}}
\def\ean{\end{align*}}

\def\be{\begin{eqnarray}}
\def\ee{\end{eqnarray}}
\def\ben{\begin{eqnarray*}}
\def\een{\end{eqnarray*}}

\def\bqq{\begin{equation}}
\def\eqq{\end{equation}}
\def\bqqn{\begin{equation*}}
\def\eqqn{\end{equation*}}

\def\elabel#1{\label{e:#1}}

\def\sq{$\Box$}

\def\qed{\ifmmode\sq\else{\unskip\nobreak\hfil
\penalty50\hskip1em\null\nobreak\hfil\sq
\parfillskip=0pt\finalhyphendemerits=0\endgraf}\fi\par\medbreak}

\newsavebox{\junk}
\savebox{\junk}[1.6mm]{\hbox{$|\!|\!|$}}

\def\til={{\widetilde =}}

 \def\eq#1/{(\ref{#1})}

\newtheorem{theorem}{Theorem}[section]

\newtheorem{proposition}[theorem]{Proposition}
\newtheorem{lemma}[theorem]{Lemma}

\def\eq#1/{(\ref{e:#1})}

\newcommand{\beqn}[1]{\notes{#1}%
\begin{eqnarray} \elabel{#1}}

\newcommand{\eeqn}{\end{eqnarray} } 

\newcommand{\beq}[1]{\notes{#1}%
\begin{equation}\elabel{#1}}

\newcommand{\eeq}{\end{equation}} 

\def\bdes{\begin{description}}
\def\edes{\end{description}}

\def\notes#1{}

\definecolor{mag}{rgb}{0.7,0,0.3}
\definecolor{dgreen}{rgb}{0.1,0.5,0.1}
\definecolor{dred}{rgb}{.8,0,0}
\definecolor{gray}{rgb}{.8,.8,.8}
\definecolor{brown}{rgb}{0.6451,0.3706,0.1745}

\begin{document}
 
\title{\vspace{-1.0cm}%
Entropy and the discrete central limit theorem}

\author
{
	Lampros Gavalakis
    \thanks{Department of Engineering,
	University of Cambridge,
        Trumpington Street,
	Cambridge CB2 1PZ, U.K.
                Email: \texttt{\href{mailto:lg560@cam.ac.uk}%
			{lg560@cam.ac.uk}}.
	L.G.\ was supported in part by EPSRC grant number RG94782.
        }
\and
        Ioannis Kontoyiannis 
    \thanks{Statistical Laboratory, DPMMS,
	University of Cambridge,
	Centre for Mathematical Sciences,
        Wilberforce Road,
	Cambridge CB3 0WB, U.K.
                Email: \texttt{\href{mailto:yiannis@maths.cam.ac.uk}%
			{yiannis@maths.cam.ac.uk}}.
	I.K.\ was supported in part by the Hellenic Foundation for Research 
	and Innovation (H.F.R.I.) under the ``First Call for H.F.R.I. Research 
	Projects to support Faculty members and Researchers and the 
	procurement of high-cost research equipment grant,'' project 
	number 1034.
        }
}

\date{\today}

\maketitle

\begin{abstract}
A strengthened version of the central limit theorem for 
discrete random variables is established, relying only 
on information-theoretic tools and elementary arguments.
It is shown that the relative entropy between the standardised 
sum of $n$ independent and identically distributed 
lattice random variables and an appropriately discretised Gaussian,
vanishes as $n\to\infty$. 
\end{abstract}


\medskip



\section{Introduction}

Suppose $X_1,X_2,\ldots$ are zero-mean,
independent and identically distributed
(i.i.d.), continuous random variables, with finite variance
$\sigma^2$. The study of the entropy $h(\hat{S}_n)$ of the
standardised sums $\hat{S}_n=\frac{1}{\sqrt{n}}\sum_{i=1}^n X_i$ 
has a long history, starting with the 1959 work of
Linnik~\cite{linnik:59}.
Recall that the {\em entropy} of a continuous
random variable $Y$ with density $f$ 
is $h(Y)=-\int f\log f$, where `$\log$' denotes 
the natural logarithm.

Barron~\cite{barron:clt} 
showed that, as $n\to\infty$,
\be
h(\hat{S}_n)\to
h(Z)=\frac{1}{2}\log(2\pi e\sigma^2),
\label{eq:barron}
\ee
where $Z\sim N(0,\sigma^2)$
is a zero-mean Gaussian
random variable
with variance~$\sigma^2$.
Barron's proof combined
earlier results by 
Brown~\cite{brown:82}
together with an integral form 
of de Bruijn's identity for the entropy
and a convolution 
inequality for the Fisher information
\cite{stam:59,blachman:65}.

The fact that the Gaussian
has maximal entropy among all random variables with variance no greater
than~$\sigma^2$ invites an appealing analogy between~(\ref{eq:barron})
and the second law of thermodynamics.
Indeed, this analogy
was carried further when it was shown that
the entropy $h(\hat{S}_n)$ in fact {\em increases}
to the maximum entropy $h(Z)=\frac{1}{2}\log(2\pi e\sigma^2)$.
This was first established using analytical tools 
by Artstein~et~al.~\cite{artstein:04},
and later using information-theoretic techniques
in~\cite{tulino:06} and~\cite{madiman:07}.

Let $D(f\|g)=\int f\log(f/g)$
denote the
{\em relative entropy} between two probability densities
$f,g$ on $\RL$.
For a continuous random variable $Y$ with density $f$
and variance $\sigma^2<\infty$,
we write $D(Y)$ for $D(f\|\phi)$,
the relative
entropy between $f$ and the Gaussian density $\phi$ 
with the same mean and variance as $Y$. 
Then we always have,
\be
D(Y)=\frac{1}{2}\log(2\pi e\sigma^2)-h(Y),
\label{eq:expand}
\ee
which implies that the convergence
of $h(\hat{S}_n)$ to $h(Z)$ is equivalent
to,
\be
D(\hat{S}_n)\to 0,\qquad\mbox{as}\;n\to\infty.
\label{eq:artstein}
\ee
In view of Pinsker's inequality,
$2\|f-g\|_{\rm TV}^2\leq D(f\|g)$~\cite{csiszar:67,kullback:67},
the relative entropy convergence in
(\ref{eq:artstein})~is enough, e.g., to recover the central limit theorem 
(CLT)
in the sense of total variation convergence.

Note that, not only do the results 
in~(\ref{eq:barron}) and~(\ref{eq:artstein})
not rely on the CLT, but they
imply a strong form of the 
CLT, established without using 
any of 
the usual probabilistic techniques.

In the case of 
{\em discrete} random variables $\{X_n\}$,
there is no immediately obvious starting point 
for identifying a corresponding
connection between the CLT and the 
entropy of the standardised sums $\hat{S}_n$;
for example,
the distribution of $\hat{S}_n$ is 
orthogonal to the Gaussian and the 
relative entropy between them is 
always infinite. 
The main contribution of this work
is the development of natural discrete 
analogs of the ``entropic'' CLTs
in~(\ref{eq:barron}) and~(\ref{eq:artstein}).

For i.i.d.\ random variables $\{X_n\}$,
write $S_n$ for the partial sums $X_1+X_2+\cdots+X_n$,
so that $\hat{S_n}=\frac{1}{\sqrt{n}}S_n.$
If the $\{X_n\}$ are continuous with finite variance, 
then by the elementary scaling property of 
the entropy~\cite{cover:book2},~(\ref{eq:barron})
can equivalently be written,
\be
\lim_{n\to\infty}\Big[h(S_n)-\log \sqrt{n}\Big]
=\frac{1}{2}\log(2\pi e\sigma^2).
\label{eq:entropycC}
\ee
The {\em entropy} of 
a discrete random variable $Y$ with probability 
mass function $p$ on a countable set $A$ is
$H(Y)=-\sum_{y\in A}p(y)\log p(y)$. 
Our first result is the analog 
of~(\ref{eq:entropycC}) for lattice random variables.
We say that $Y$ has a lattice distribution with span
$h>0$ if its support is a subset of $\{a+kh\;:\;k\in\IN\}$
for some $a\in\RL$;
the span $h$ is {\em maximal} if it is the largest such $h$.

\begin{theorem}[Entropy convergence]
\label{thm:entropycD}
If $S_n$, $n\geq 1$, are the partial sums
of a sequence $\{X_n\}$ of i.i.d.\ lattice random variables
with finite variance $\sigma^2$ and maximal span $h$,
then:
\be
\lim_{n\to\infty}\Big[H(S_n)-\log \frac{\sqrt{n}}{h}\Big]
=\frac{1}{2}\log(2\pi e\sigma^2).
\label{eq:entropycD}
\ee
\end{theorem}

Since the discrete entropy does not scale in the same way
as the continuous entropy (e.g.,
$H(S_n)=H(\hat{S}_n)$), 
the equivalence between the convergence in~(\ref{eq:entropycD}) 
and a discrete version of the entropic CLT
$D(\hat{S}_n)\to 0$ 
is no longer immediate.
Nevertheless,
it is possible to establish a result analogous
to that in~(\ref{eq:artstein})
in the discrete case, 
as shown in Theorem~\ref{thm:relentropycD} below.

For discrete random variables $X,Y$ with probability
mass functions $p,q$, respectively,
on the same countable set $A$, the {\em relative entropy}
$D(p\|q)$ 
between $p$ and $q$ 
is defined as $D(p\|q)=\sum_{x\in s}p(x)\log(p(x)/q(x))$,
where the sum is interpreted as the Lebesgue integral 
of $\log(p/q)$ with
respect to the probability measure induced by $p$ on $\RL$.

\newpage

Suppose $Y$ is lattice random variable
maximal span~$h$,
values in $A=\{a+kh\;:\;k\in\IN\}$,
mean $\mu$, 
and finite variance $\sigma^2$.
We write $D(Y)$ for the relative entropy $D(p\|q)$
between 
the probability mass function $p$ of $Y$
and the probability mass function $q$
of a Gaussian random variable 
$Z\sim N(\mu,\sigma^2)$ 
quantised on $A$ as,
\be
q(a+kh)=\int_{a+kh}^{a+(k+1)h}\phi(x)dx,\qquad k\in\IN,
\label{eq:quantise}
\ee
where $\phi$ is the $N(\mu,\sigma^2)$ density.
Observe that, by definition, $D(Y+c)=D(Y)$ for any constant~$c$.

\begin{theorem}[Discrete entropic CLT]
\label{thm:relentropycD}
If $\hat{S}_n$, $n\geq 1$, are the standardised sums
of a sequence $\{X_n\}$ of i.i.d.\ lattice random variables
with finite variance, then:
\be
D\Big(\frac{1}{\sqrt{n}}\sum_{i=1}^n(X_i-\mu)\Big)
=D(\hat{S}_n)\to 0,\qquad\mbox{as}\;n\to\infty.
\label{eq:relentropycD}
\ee
\end{theorem}

As in the continuous case, Pinsker's inequality combined
with the triangle inequality for the total variation norm
imply a strong version of the CLT: 
Taking $\mu=0$ without loss of generality,
let $Z\sim N(0,\sigma^2)$
and let $Z_n$ be the quantised Gaussian as in the
definition of $D(\hat{S}_n)$. Then, 
$$\|\hat{S}_n-Z\|_{\rm TV}\leq
\sqrt{\frac{1}{2}D(\hat{S}_n)}+\|Z_n-Z\|_{\rm TV}\to 0,
\qquad
\mbox{as}\;n\to\infty,$$
since the first term vanishes by Theorem~\ref{thm:relentropycD}
and the second term vanishes by the definition of~$Z_n$.
Alternatively, the fact that
$\|\hat{S}_n-Z_n\|_{\rm TV}\to 0$ 
implied by Theorem~\ref{thm:relentropycD}
readily translates to local-CLT-like results.

\medskip

\noindent
{\bf Paper outline and proof ideas.}
In the end of this Introduction we discuss the intriguing
connection between the CLT
and Shannon's entropy power inequality.
In Section~\ref{prelimsection} we prove
that the entropic CLT statements~(\ref{eq:entropycD})
and~(\ref{eq:relentropycD}) in
Theorems~\ref{thm:entropycD} 
and~\ref{thm:relentropycD},
respectively, are equivalent.
There we also establish a 
relation between
$D(\hat{S}_n)$ and $D(\hat{S}_n+U)$, when 
$\hat{S}_n$ are the standardised sums of lattice
random variables $X_i$ and and $U$ is an
appropriate (continuous) uniform random variable,
independent of the $X_i$
(Lemma~\ref{discretecontinuousD}).
In Section \ref{demoivresection} we 
establish two special cases of Theorem~\ref{thm:entropycD};
first, when the $X_i$ in $S_n=\sum_{i=1}^nX_i$ 
are symmetric Bernoulli random variables,
$X_i\sim{\rm Bern}(1/2)$, and then when 
each $X_i$ can be written 
as the sum $X_i=V_i+B_i$
of a lattice random variable~$V_i$
with maximal span $h=1$ and a
$B_i\sim{\rm Bern}(1/2)$
independent of~$V_i$.

At first sight it might be tempting to hope
that the result of Theorem~\ref{thm:relentropycD} 
could be derived 
from its
continuous counterpart~(\ref{eq:artstein})
via a simple quantisation argument using the 
``data processing'' property of
relative entropy~\cite{cover:book2}, but this does not appear
to be the case.  Instead, Barron's
continuous result~(\ref{eq:artstein}) 
is employed in a more indirect way in the proof 
of the special case of Theorem~\ref{thm:entropycD}
given in Section~\ref{demoivresection}.
This is then used in the proof of our main
result, the general case of Theorem~\ref{thm:entropycD},
in Section~\ref{mainsection}.
In addition to Barron's result~(\ref{eq:artstein}),
and to simple information-theoretic properties 
and some well-known bounds and identities for 
the Fisher information, the other main ingredient in 
the proof of Theorem~\ref{thm:entropycD} is
an elementary technique known as 
``Bernoulli part decomposition,''
described in Section~\ref{mainsection}.
Theorem~\ref{thm:relentropycD} is an
immediate consequence of Theorem~\ref{thm:entropycD}
combined with Theorem~\ref{discreteDH}.

\medskip

The CLT for discrete random variables has been investigated
from an information-theoretic point of view
by, among others, Shimizu \cite{shimizu:75} and Brown~\cite{brown:82},
who obtained the convergence of $\hat{S}_n$ 
to a Gaussian in distribution 
(but not for entropy or relative entropy) by 
proving convergence for the Fisher information of smoothed 
versions of $\hat{S}_n$.
The Bernoulli part decomposition technique was first used (implicitly)
by Mineka~\cite{mineka:73} in a different context,
and by McDonald~\cite{mcdonald:80} 
and Davis and McDonald~\cite{davis:95},
who derived conditions under which the standardised sums 
of independent discrete random variables
satisfy the local CLT~\cite{esseen:45,gnedenko:48}.
In the reverse direction, Takano~\cite{takano:87} used
the local CLT to derive entropy expansions for the
standardised sums $\hat{S}_n$ as in our Theorem~\ref{thm:entropycD}.

There is a significant line of work re-examining
core probabilistic results through the lens of 
information theory.
In terms of ideas as well as techniques,
perhaps the works closest in spirit to the present development
are those providing information-theoretic treatments of
Poisson approximation 
\cite{harremoes:01,konto-H-J:05}
and compound Poisson approximation \cite{Kcompound:10,jKm:13}.

\medskip

\noindent
{\bf The CLT and the entropy power inequality.}
The earliest indication of a nontrivial connection between
the CLT and information-theoretic ideas comes from
Shannon's entropy power 
inequality (EPI)~\cite{shannon:48,stam:59,blachman:65}. 
For i.i.d.\ continuous random variables $X_1,X_2$,
the EPI states that,
$$h(X_1+X_2)\geq h(X_1)+\frac{1}{2}\log 2,$$
with equality if and only if $X_1,X_2$ are Gaussian.
Using the scaling property of the
entropy, this implies that $h(\hat{S}_{2n})\geq h(\hat{S}_n)$
for all $n$,
which likely provided some of the initial motivation for 
the works \cite{linnik:59,shimizu:75,brown:82,barron:clt} mentioned earlier. 
Further, a generalisation of the EPI was used 
to prove the monotonic increase of the entire
sequence $\{h(\hat{S}_n)\}$ to $\frac{1}{2}\log (2\pi e\sigma^2)$
in~\cite{madiman:07}.

For i.i.d.\ discrete random variables $X_1,X_2$, it is easy to 
see (by considering random variables with entropy 
close to zero) that the obvious discrete analog,
$H(X_1 + X_2) \geq H(X_1) + \frac{1}{2}\log{2},$ fails to hold
in general. On the other hand, 
Tao~\cite{tao:10} showed that, 
for any $\epsilon>0$,
\begin{equation*} 
H(X_1 + X_2) \geq H(X_1) + \frac{1}{2}\log{2} - \epsilon,
\end{equation*} 
for all i.i.d.\ pairs $X_1,X_2$ such that $H(X_1)$ is
large enough depending on $\epsilon$. Tao's proof relies on the 
inverse sumset theory for entropy 
developed in~\cite{tao:10}.
A careful examination of the proof shows that 
not only is the lower bound on $H(X_1)$ 
at least,
\vspace*{-0.05in}
\be
\Omega\Bigl(\frac{1}{\epsilon}^{\frac{1}{\epsilon}^%
{\frac{1}{\epsilon}}}\Bigr),
\label{eq:Omega}
\ee
but the implied absolute constants 
are also very large. 
Although the 
answer to the natural
question of how much this bound can be
improved remains unclear
(see, e.g.,~\cite{abbe:14}
for some related bounds),
in view of 
the results in this paper,
particularly
the nonasymptotic versions
of Theorem~\ref{discreteDH} and Lemma~\ref{discretecontinuousD},
we expect that perhaps
if one restricts attention to 
lattice random variables with finite variance, 
it may be possible to significantly improve on~(\ref{eq:Omega}).

Interestingly,
Tao~\cite{tao:10} further conjectured that,
for any $n\geq 2$, 
\begin{equation} \label{taosconjecture}
H(X_1+\ldots+X_n) \geq H(X_1+\ldots+X_{n-1}) 
+  \frac{1}{2}\log\Big(\frac{n}{n-1}\Big) - \epsilon,
\end{equation} 
as long as $H(X_1)$ is sufficiently large depending on $n$ 
and $\epsilon$. 
The present results again suggest that
this conjecture might be easier 
to prove if attention is restricted to lattice
random variables with finite variance.
Specifically, 
in this setting~(\ref{taosconjecture})
can be interpreted as an ``approximate monotonicity'' 
refinement of our Theorems~\ref{thm:entropycD}
and~\ref{thm:relentropycD}:
By the nonasymptotic form of Theorem~\ref{discreteDH} 
and the fact that $H(X_1) \to \infty$ 
implies $\VAR(X_1) \to \infty$, for lattice $X_1$ with finite 
variance the conjecture~\eqref{taosconjecture} is equivalent to:
$$
D(\hat{S}_n)\leq D(\hat{S}_{n-1})+\epsilon.
$$

\section{Entropy, relative entropy, and Fisher information}
\label{prelimsection}

Let $X_1,X_2,\ldots,$ be i.i.d.\  lattice random variables with 
values in $\{a + kh: k \in \mathbb{Z}\}$,
mean $\mu$,
and finite variance $\sigma^2.$
As before, write $S_n$ for the partial sums
$\sum_{i=1}^n{X_i}$ and $\hat{S}_n$ for the standardised
sums $\frac{1}{\sqrt{n}}S_n$,
and recall the definition of the 
relative entropy $D(Y)$ between a lattice
random variable $Y$ and and an appropriately
quantised Gaussian as in~(\ref{eq:quantise}).

Our first observation is that the ``entropy deficit,''
$$\frac{1}{2}\log{(2\pi e\sigma^2)} 
- \Bigl[H(S_n) - \log{\frac{\sqrt{n}}{h}}\Bigr],$$
can be viewed as a measure of the ``Gaussianity'' of 
the lattice sum $S_n$. 
Theorem~\ref{discreteDH}
shows that the
entropic CLTs stated in Theorems~\ref{thm:entropycD}
and~\ref{thm:relentropycD} are equivalent.

\begin{theorem}[Entropy and relative entropy solidarity]
\label{discreteDH} \;
Suppose $\{X_n\}$ are i.i.d.\ lattice random variables 
with finite variance $\sigma^2$ and maximal
span $h>0$.
Then the partial sums $S_n$ and the
standardised sums $\hat{S}_n$ of the $X_i$
satisfy, as $n\to\infty$,
$$
D(\hat{S}_n)
= \frac{1}{2}\log{(2\pi e\sigma^2)} - 
\Bigl[H(S_n) - \log{\frac{\sqrt{n}}{h}}\Bigr]   
+ O\Bigl(\frac{1}{\sqrt{n}}\Bigr).
$$
In fact, for all $n\geq 1$, the $O(1/\sqrt{n})$ error term
is absolutely bounded by:
$$\frac{h}{\sigma \sqrt{n}}\Big[1 + \frac{h}{2 \sigma\sqrt{n}}\Big].$$
\end{theorem}

\noindent
{\sc Proof. } 
Because $H(Y)$ is translation invariant and, as
noted in the introduction, so is $D(Y)$, we may 
assume that $\mu=0$ without loss of generality.

Since the $X_i$ take values in $\{a+kh\;:\;k\in\IN\}$,
$S_n$ takes values in $\{na+kh\;:\;k\in\IN\}$
and $\hat{S}_n$ in $\{\sqrt{n}a+kh/\sqrt{n}\;:\;k\in\IN\}$.
Let $p,q$ denote the probability mass
functions of $S_n$ and of the quantised
Gaussian in the definition of $D(\hat{S}_n)$,
respectively. Writing $\phi$ for the 
standard normal density, 
for each $k\in\IN$ we have,
$$
q\Big(\sqrt{n}a+k\frac{h}{\sqrt{n}}\Big)=
\int_{\sqrt{n}a+kh/\sqrt{n}}^{\sqrt{n}a+(k+1)h/\sqrt{n}}
\frac{1}{\sigma}\phi\Big(\frac{x}{\sigma}\Big)dx
=
\frac{h}{\sigma\sqrt{n}}
\phi\Big(\frac{n a+\xi_k h}{\sigma\sqrt{n}}\Big),
$$
for some $\xi_k\in[k,k+1]$.
Using this, we can bound the absolute difference,
\begin{align*} 
\Delta_n
&:=
	\left| 
	D(\hat{S}_n)-\Big[
	\frac{1}{2}\log{(2\pi e\sigma^2)} - 
	H(S_n) + \log{\frac{\sqrt{n}}{h}}
	\Bigr] 
	\right| \\
&
	\;=\left|
	\sum_{k \in \mathbb{Z}}p(na + kh)
	\log \left(
	\frac{p(na + kh)}{\frac{h}{\sigma\sqrt{n}}
		\phi(\frac{na + \xi_kh }{\sigma\sqrt{n}})}\right)
	- \sum_{k \in \mathbb{Z}}
	p(na + kh)
	\log\left(\frac{p(na + kh)}
	{\frac{h}{\sigma \sqrt{n}}\phi(\frac{na + kh}
	{\sigma\sqrt{n}})}\right)
	\right|,
\end{align*}
where the second sum contains the last three terms in $\Delta_n$.
Simplifying we obtain,
\begin{align*}
\Delta_n
&\leq \frac{1}{2n\sigma^2}\sum_{k\in\IN}
{p(na + kh) \Bigl|(na + kh)^2-(na + \xi_kh)^2\Bigr|} \\
&\leq \frac{1}{2n\sigma^2}\sum_{k\in\IN}
{p(na + kh)\bigl(2h|na + kh | + h^2\bigr)} 
\;\leq\; \frac{h}{\sigma \sqrt{n}} + \frac{h^2}{2 \sigma^2 n},
\end{align*}
as required, where the last step follows from the 
Cauchy-Schwarz inequality 
and the fact that the variance of $S_n$ is $n\sigma^2$.
\qed

By the nonnegativity of relative entropy we obtain
the following standard upper bound, which can be viewed
as a discrete analog of the maximum entropy property 
of the Gaussian:
\begin{equation} 
\label{entropyupperbound}
\Bigl[H(S_n) - \log{\frac{\sqrt{n}}{h}}\Bigr]   
\leq \frac{1}{2}\log{(2\pi e\sigma^2)} 
+ 
\frac{h}{\sigma \sqrt{n}}\Big[1 + \frac{h}{2 \sigma\sqrt{n}}\Big].
\end{equation}
In fact, we
can easily obtain a stronger bound.
Let $U$ be
an independent uniform random variable 
on $(-1/2,1/2)$. Then,
by the definitions of the continuous 
and discrete entropies,
\begin{equation} 
\label{eq:uniform}
H(S_n) - \log{\frac{\sqrt{n}}{h}} 
= h\Big(\hat{S}_n + \frac{h}{\sqrt{n}}U\Big).
\end{equation}
And using the maximum 
maximum entropy
property
of the Gaussian yields:

\begin{proposition}
If $S_n$ is the sum of $n$ i.i.d.\ lattice random variables
with maximal span $h>0$ and finite variance $\sigma^2$,
then:
\begin{equation} 
\label{entropyupperboundbyuni}
H(S_n) - \log{\frac{\sqrt{n}}{h}} 
\leq \frac{1}{2}\log\Big[2\pi e \Big(\sigma^2 + \frac{h^2}{12n}\Big)\Big].
\end{equation}
\end{proposition}

In the special case $h=1$, $n=1$,
the bound~(\ref{entropyupperboundbyuni}) appeared 
in~\cite{massey:88}.
It was recently exploited 
further in~\cite{rioul:21}, where an improved 
inequality was also established for large $\sigma^2$
via the Poisson summation formula. 
For any $n$ and $h=1$, 
\eqref{entropyupperboundbyuni} also appeared in \cite{bobkov:20},
as a special case of an inequality for R\'enyi entropies.

The following lemma will be used in the proof 
of Theorem~\ref{bernsmooth}.
It highlights the asymptotic equivalence between 
the discrete and continuous versions of the relative
entropy $D(\hat{S}_n)$.

\begin{lemma} 
\label{discretecontinuousD}
Under the assumptions of Theorem~\ref{discreteDH}, 
let $U$ be an independent uniform random variable 
on $(-1/2,1/2).$ 
Then, as $n \to \infty$,
\begin{equation*} 
D(\hat{S}_n) 
= D\Big(\hat{S}_n + \frac{h}{\sqrt{n}}U\Big) 
+ O\Big(\frac{1}{\sqrt{n}}\Big).
\end{equation*}
In fact, for all $n\geq 1$, the $O(1/\sqrt{n})$
error term is absolutely bounded by:
$$
\frac{h}{\sigma\sqrt{n}}\Big[1 + \frac{13h}{24\sigma\sqrt{n}}\Big].
$$
\end{lemma}

\noindent
{\sc Proof. }
As in the proof of Theorem~\ref{discreteDH},
we may assume without loss of generality
that the $X_i$ have zero mean.
Note that $\hat{S}_n+(h/\sqrt{n})U$ has variance
$\sigma^2+h^2/(12n)$.
Using the finite-$n$ bound in Theorem~\ref{discreteDH}
and the general property~(\ref{eq:expand}) of the
relative entropy,
\begin{align*}
&
	\left|
	D(\hat{S}_n) 
	-D\Big(\hat{S}_n + \frac{h}{\sqrt{n}}U\Big)
	\right|\\
&
	\leq 
	\left|
	D\Big(\hat{S}_n + \frac{h}{\sqrt{n}}U\Big) 
	- \frac{1}{2}\log{(2\pi e\sigma^2)} 
	+ \Bigl[H(S_n) - \log{\frac{\sqrt{n}}{h}}\Bigr]\right|   
	+ \frac{h}{\sigma\sqrt{n}} + \frac{h^2}{2 \sigma^2 n}\\
&
	=
	\left|
	\frac{1}{2}\log\Big(1+
	\frac{h^2}{12n\sigma^2}\Big)
	-
	h\Big(\hat{S}_n + \frac{h}{\sqrt{n}}U\Big) 
	+ \Bigl[H(S_n) - \log{\frac{\sqrt{n}}{h}}\Bigr]\right|   
	+ \frac{h}{\sigma\sqrt{n}} + \frac{h^2}{2 \sigma^2 n},
\end{align*}
and using~(\ref{eq:uniform}),
$$	\left|
	D(\hat{S}_n) 
	-D\Big(\hat{S}_n + \frac{h}{\sqrt{n}}U\Big)
	\right|
\leq
	\frac{1}{2}\log\Big(1+
	\frac{h^2}{12n\sigma^2}\Big)
+ \frac{h}{\sigma\sqrt{n}} + \frac{h^2}{2 \sigma^2 n}
\leq  \frac{h}{\sigma\sqrt{n}} + \frac{13h^2}{24\sigma^2n},
$$
where the last inequality follows from
the elementary bound $\log{(1+x)} \leq x$, $x > 0.$ 
\qed

We close this section by recalling some simple
convolution inequalities that will be 
used in the following sections.
If $X,Y$ are independent discrete random variables,
then
\cite{cover:book2}:
\begin{equation}\label{entropyadd}
H(X+Y) \geq H(X).
\end{equation}
Similarly, if $X$ is a continuous random variable 
and $Y$ an arbitrary independent random variable,
then
\cite{cover:book2}:
\begin{equation}\label{differentialentropyadd}
h(X+Y) \geq h(X).
\end{equation}
Finally, for a continuous random variable $X$
with a continuously differentiable density $f$,
we define the {\em Fisher information} of $X$ as
$I(X) = \int (f')^2/f$.
If the independent random variables $X,Y$ 
have continuously differentiable densities with bounded derivatives,
then~\cite[Lemma 5.5]{brown:82}:
\begin{equation}
\label{fisheradd}
I(V+W) \leq I(V).
\end{equation}

\section{Binomial sums and Bernoulli smoothing}
\label{demoivresection}

We first establish
a nonasymptotic version of 
Theorem~\ref{thm:entropycD} in the 
special case when $S_n$
is the sum
of independent ${\rm Bern}(1/2)$ random variables,
so that 
$S_n\sim{\rm Bin}(n,1/2)$
has
a binomial distribution with parameters $n$ and $1/2$. 
Although this elementary result 
is largely known~\cite{chang:79,hughes:96,jacquet:99},
we state and prove it explicitly as it is the first 
step towards the proof of our main result,
Theorem~\ref{thm:entropycD}. Also, as earlier proofs 
of~(\ref{binomialentropy}) actually
use the CLT, relying on such arguments would defeat 
our main claim,
namely, that of obtaining a complete proof of the 
entropic CLT without using any of the standard probabilistic
normal approximation techniques 
or, of course, the CLT itself.

\begin{proposition}[Binomial entropy]
\label{demoivreentropy}
If $S_n\sim{\rm Bin}(n,1/2)$, then
for all $n\geq 2$:
\begin{equation} 
\left|\Bigl[H(S_n) - \log{\sqrt{n}}\Bigr]
-\frac{1}{2}\log{\Bigl(\frac{1}{2}\pi e\Bigr)}
\right|\leq\frac{4}{\sqrt{n}}.
\label{binomialentropy}
\end{equation}
\end{proposition}

Note that~(\ref{binomialentropy})
combined with the finite-$n$ version
of Theorem~\ref{discreteDH}
also  yields:
\begin{equation*} 
\label{binomialrelativeentropy}
D(\hat{S}_n) 
\leq \frac{8}{\sqrt{n}},
\qquad n\geq 2.
\end{equation*}

\noindent
{\sc Proof. }
The general upper bound in~(\ref{entropyupperboundbyuni})
in this particular case gives,
\be
H(S_n)-\log\sqrt{n}
\leq
\frac{1}{2}\log\Big(\frac{1}{2}\pi e\Big)
+\frac{1}{2}\log\Big(1+\frac{1}{12n}\Big)
\leq
\frac{1}{2}\log\Big(\frac{1}{2}\pi e\Big)+\frac{1}{24n}.
\label{eq:binup}
\ee
For the proof of the corresponding
lower bound we only consider even $n$; the case
of odd $n$ is similar.
Let 
$b_n(k) = {n\choose k}2^{-n}$
denote the ${\rm Bin}(n,1/2)$ probabilities,
and for fixed $n\geq 2$ write $a_k = b_n(\frac{n}{2}+k)$,
for $-n/2\leq k\leq n/2$.
Following a simple argument by Feller~\cite[VII, 2]{fellerI:book}, 
we first observe that for $k\geq 1$,
$$a_k = a_0\times\frac{\frac{n}{2}(\frac{n}{2}-1)\cdots
(\frac{n}{2}-k+1)}{(\frac{n}{2}+1)(\frac{n}{2} + 2)\cdots(\frac{n}{2} + k)} 
= a_0\times\frac{(1-\frac{2}{n})(1-\frac{4}{n})\cdots
(1-\frac{2(k-1)}{n})}{(1+\frac{2}{n})\cdots(1+\frac{2k}{n})},
$$
and then use the elementary bounds,
$1-x \leq e^{-x}$ 
and $1 + x \geq e^{x-x^2}$, for $x \in [0,1),$ 
to obtain that,
\begin{equation*}
a_k \leq a_0
\exp\Big\{
-\frac{4}{n}\bigl[1+\cdots+(k-1)\bigr] - \frac{2k}{n} + 3\frac{k^3}{n^2}
\Big\} 
= a_0e^{-\frac{2k^2}{n} + \frac{3k^3}{n^2}}.
\end{equation*}
By Robbins' finite-$n$ version of Stirling's formula, 
e.g.~\cite[II, (9.15)]{fellerI:book}, we can easily bound,
$$
a_0 \leq \Big(\frac{\pi n}{2}\Big)^{-1/2}e^{\frac{1}{12n}},
$$
so that, for $k\geq 0$,
\begin{equation}
a_k \leq \Big(\frac{\pi n}{2}\Big)^{-1/2}
e^{-\frac{2k^2}{n}}e^{\frac{3k^3}{n^2} 
+ \frac{1}{12n}}.
\label{eq:feller}
\end{equation}
Since $a_k=a_{-k}$, the same bound
holds for all $-n/2\leq k\leq n/2$,
with $|k|$ in place of $k$.
And substituting~(\ref{eq:feller}) into the 
logarithmic term in the definition of $H(S_n)$
gives,
\begin{align}\nonumber
H(S_n) 
&= -\sum_{k=-n/2}^{n/2}a_k\log a_k \\ \nonumber
&\geq \log{\sqrt{n}} 
	+\frac{1}{2}\log\Big(\frac{1}{2}\pi\Big)
	+ \frac{2}{n}\sum_{k = n/2}^{n/2}
	a_k k^2
	- \frac{3}{n^2}\sum_{k = -n/2}^{n/2}
	a_k |k|^3 - \frac{1}{12n} \\ 
&\geq 
	\log{\sqrt{n}} + \frac{1}{2}\log\Big(\frac{1}{2}\pi\Big)
	+ \frac{1}{2} - \frac{4}{\sqrt{n}},
	\label{hbnlowerbound}
\end{align}
where in the last step we used
the fact that the variance of $S_n$ is $n/4$
and its third absolute central moment
is bounded above by $n^{3/2}$.

The result follows from~(\ref{eq:binup})
and~\eqref{hbnlowerbound}.
\qed

Next, we extend the result of
Proposition~\ref{demoivreentropy}
to the case when each
$X_i$ in $S_n$
can be written as the independent
sum $X_i=V_i+B_i$ of a lattice random 
$V_i$ and a $B_i\sim{\rm Bern}(1/2)$. 
The proof of Theorem~\ref{bernsmooth}
is a key step towards 
the proof of the general case 
of Theorem~\ref{thm:entropycD} 
in the next section.
We refer to the addition of an independent
Bernoulli to a lattice random variable
as ``Bernoulli smoothing,''
in analogy to the Gaussian smoothing
step used in~\cite{shimizu:75,brown:82,barron:clt}
along the corresponding development
in the continuous case. There,
one considers $X_i + \sqrt{t}Z_i$,
where the $Z_i$ are 
standard normals, so that the
resulting random variables have
differentiable densities that smoothly interpolate
between the distribution of $X_i$ and 
the Gaussian, as $t$ varies.
In our case, the addition of a
binomial random variable to the
partial sums $S_n$ facilitates the 
use of Proposition~\ref{demoivreentropy}, and also 
allows us to establish a uniform integrability
property which can be used to 
exploit the fact that Fisher information 
decreases on convolution. 

\begin{theorem}[Bernoulli smoothing]
\label{bernsmooth}
Suppose $\{V_n\}$ are i.i.d.\ lattice random variables 
with finite variance $\sigma_V^2$ and maximal span $h=1$,
and let $\{B_n\}$ be i.i.d.\ ${\rm Bern}(1/2)$, 
independent of $\{V_n\}$.
Then:
\begin{equation*}
\lim_{n \to \infty}{\Biggl[H\Biggl(
\sum_{i=1}^n\big[V_i + B_i\big]\Biggr) - \log{\sqrt{n}}\Biggr]} 
= \frac{1}{2}\log{\Bigl(2\pi e\Big(\sigma_V^2 + \frac{1}{4}\Big)\Bigr)}.
\end{equation*}
\end{theorem}

For a continuous random variable with continuously differentiable
density $f$, the {\em score function} $\rho$ of $Y$ is 
$\rho=f'/f$, so that the Fisher information $I(Y)$ can
be expressed $I(Y)=\int f\rho^2$. In particular,
if $\phi$ is the $N(\mu,\sigma^2)$ density,
then its score function $\rho_\phi$ is linear,
$\rho_\phi(z)=-(z-\mu)/\sigma^2$, $z\in\RL$. 
For the proof of Theorem~\ref{bernsmooth},
we will find it convenient to use
the {\em standardised Fisher information} $J(Y)$,
which, when $Y$ has mean $\mu$ and variance $\sigma^2$,
is defined as $J(Y)=\sigma^2\int f(\rho-\rho_\phi)^2$,
or, equivalently, 
\be
J(Y) = \sigma^2I(Y) - 1.
\label{eq:SFI}
\ee


\newpage

\noindent
{\sc Proof. }
Let $U$ be an independent random variable,
uniformly distributed on $(-1/2,1/2)$,
and write $S_n$ for the binomial sum
$S_n=\sum_{i=1}^nB_i$. In view of
Theorem~\ref{discreteDH} and Lemma~\ref{discretecontinuousD}, 
it suffices to show that, as $n\to\infty$,
$$D\Bigl(\frac{1}{\sqrt{n}}
\sum_{i=1}^n{V_i} + \frac{1}{\sqrt{n}}S_n + \frac{1}{\sqrt{n}}U\Bigr)
\rightarrow 0.$$
Using Barron's integral form of de Bruijn's 
identity~\cite[Eq. (4.1)]{barron:clt}, this can be
expressed as,
\begin{align}
D\left(\frac{1}{\sqrt{n}}\left[\sum_{i=1}^n{V_i} 
+ S_n + U\right]\right)
=&
	\;D\left(\frac{1}{\sqrt{2n}}\left[\sum_{i=1}^n{V_i} 
	+ S_n + U\right] + \frac{1}{\sqrt{2}}Z\right) 
	\label{debruijn1}\\
&
	\;+ \int_0^{1/2}{J\Biggl(\sqrt{\frac{1-t}{n}}
	\left[\sum_{i=1}^n{V_i} 
	+ S_n + U\right] + \sqrt{t}Z\Biggr)\frac{dt}{2(1-t)}},
	\label{debruijn2}
\end{align}
where $Z$ is an independent normal random variable with the same 
mean and variance as,
$$\frac{1}{\sqrt{n}}\sum_{i=1}^n{V_i} 
+ \frac{1}{\sqrt{n}}S_n + \frac{1}{\sqrt{n}}U.$$
Writing $\mu_V,\sigma_V^2$ for the mean and variance of the
$V_i$, respectively, $Z$ can be expressed,
$Z=\frac{1}{\sqrt{n}}\sum_{i=1}^nZ_i+W_n,$
where the $Z_i$ are i.i.d.\ $N(\mu_V+1/2,\sigma_V^2+1/4)$
and $W_n\sim N(0,\frac{1}{12n})$
is independent of the $V_i$.
Therefore,
the argument of the relative entropy 
in the right-hand side of~\eqref{debruijn1} can be written,
$$
T_n=\frac{1}{\sqrt{2n}} \sum_{i=1}^n[V_i +B_i+Z_i]
+
\frac{1}{\sqrt{2n}} U
+ \frac{1}{\sqrt{2}}W_n.$$
Write
$Y_i$ for the continuous i.i.d.\ random variables
$Y_i=(V_i +B_i+Z_i)/\sqrt{2}$
and let
$\sigma_{T}^2$ and $\sigma_Y^2$
denote the variances of $T_n$
and of $Y_i$, respectively.
By~(\ref{eq:expand})
and the convolution
inequality~\eqref{differentialentropyadd}
we have, as $n\to\infty$,
\ben
D(T_n)
&=&
	\frac{1}{2}\log(2\pi e\sigma_T^2)-h(T_n)\\
&\leq&
	\frac{1}{2}\log\Big(2\pi e\Big(\sigma_Y^2+\frac{1}{12n}\Big)\Big)
	-h\Big(\frac{1}{\sqrt{n}}\sum_{i=1}^nY_i\Big)\\
&=&
	D\Big(\frac{1}{\sqrt{n}}\sum_{i=1}^nY_i\Big)+o(1),
\een
where the last relative entropy is also $o(1)$ by
the continuous entropic CLT~(\ref{eq:artstein}).
Therefore, the relative entropy in~\eqref{debruijn1} 
vanishes as $n\to\infty$, and now it suffices
to show that so does the
integral in~\eqref{debruijn2}.

An analogous argument to the one
used above for the relative entropy can be used
to show that, for each $t$, 
the standardised Fisher information
in the integrand in~(\ref{debruijn2})
vanishes with $n$. For fixed $t\in(0,1)$,
let $R_n=R_n(t)$ denote the argument of the 
standardised Fisher information
in~(\ref{debruijn2}), so that $R_n$ can
be written,
$$
R_n=\sqrt{\frac{1-t}{n}} \sum_{i=1}^n[V_i +B_i] + \sqrt{t}\hat{Z}_n
+
\sqrt{\frac{1-t}{n}} U
+ \sqrt{t}W_n,$$
where now $\hat{Z}_n\sim N(\sqrt{n}(\mu_V+1/2),\sigma_V^2+1/4).$ 
Write $Y'_i$ for the i.i.d.\ random variables
$Y'_i:=V_i+B_i$ and let
$\sigma^2_{R}$ and $\sigma_{Y'}^2$ denote the variances
of $R_n$ and $Y'_i$, respectively.
By the representation~(\ref{eq:SFI})
and the convolution 
inequality~\eqref{fisheradd},
we have that,
$$
0\leq J(R_n)
=
	\sigma^2_{R}I(R_n)-1
\leq 
	\Big(\sigma^2_{Y'}+\frac{1}{12n}\Big)
	I\Bigg(\sqrt{\frac{1-t}{n}}\sum_{i=1}^nY'_i + \sqrt{t}\hat{Z}_n\Bigg)-1,$$
which vanishes as $n\to\infty$ by 
the Fisher information convergence in~\cite[Lemma 2]{barron:clt},
since $J(\cdot)$ is translation invariant.

Finally, we show that the nonnegative sequence 
$\{J(R_n(t))\;;\;n\geq 1\}$ 
is uniformly integrable with respect to the probability measure
$\nu(dt)$ proportional to $\frac{dt}{2(1-t)}$ on $(0,1/2).$ 
In fact, we will show that it is bounded above by the uniformly
integrable sequence $\{J(R'_n(t))\}$ defined next. 

Let 
$Z'\sim N(\frac{\sqrt{n}}{2},\frac{1}{4}+\frac{1}{12n})$
and 
$Z''\sim N(\sqrt{n}\mu_V,\sigma_V^2)$ be independent
random variables such that $Z=Z'+Z''$.
Then we can write,
$$R_n=R_n'
+\sqrt{\frac{1-t}{n}}\sum_{i=1}^n V_i+\sqrt{t}Z'',$$
where,
$$R_n'=R_n'(t)=\sqrt{\frac{1-t}{n}}[S_n +U]+\sqrt{t}Z',$$
so that, by the
convolution 
inequality~\eqref{fisheradd}
and using the the representation~(\ref{eq:SFI}) twice,
$$J(R_n)=\sigma_{R}^2I(R_n)-1
\leq
(\sigma_{R'}^2+\sigma_V^2)I(R_n')-1
=\Big(1+\frac{\sigma_V^2}{\sigma_{R'}^2}\Big)J(R_n')
+\frac{\sigma_V^2}{\sigma_{R'}^2},
$$
where $\sigma_{R'}^2=\frac{1}{4}+\frac{1}{12n}$ 
is the variance of $R_n'$.

But by Proposition~\ref{demoivreentropy}, 
Lemma~\ref{discretecontinuousD} and de Bruijn's integral identity, 
$$\int_{0}^{1/2} J(R_n'(t)) \frac{dt}{2(1-t)},$$ 
vanishes as $n \rightarrow \infty.$
Therefore, 
$\{J(R'_n(t))\}$
is uniformly integrable with respect to the probability measure
$\nu(dt)\propto\frac{dt}{2(1-t)}$ on $(0,1/2)$,
and hence so is 
$\{J(R_n(t))\}$.

The result follows. 
\qed

\section{Bernoulli part decomposition}
\label{mainsection}

At the end of this section we give the proof
of Theorem~\ref{thm:entropycD}. In view of~\eqref{entropyupperbound}, 
our goal is to obtain an appropriate lower bound on the entropy
$H(S_n)$.
The main idea is to show that $S_n$ can
be asymptotically approximately decomposed 
as a sum involving a ${\rm Bin}(n,1/2)$
random variable and
then apply Theorem \ref{bernsmooth}.
The required decomposition will be based 
on the following elementary technique.

Let $X$ be an integer-valued random variable with 
probability mass function $p$ on $\IN$ and
maximal span $h=1$. The
{\em Bernoulli part decomposition} of $X$ is
the representation,
\begin{equation*}
X \equald V + WB, 
\end{equation*}
where $V$ takes values in $\IN$, $W\sim{\rm Bern}(q)$,
and $B\sim{\rm Bern}(1/2)$ is
independent of $(V,W).$
The joint probability mass function of $V$ and $W$ is given by,
\begin{align*}
p_{V,W}(k,1) &= \min\{p(k),p(k+1)\},\\
p_{V,W}(k,0) &= p(k) - \frac{1}{2}[p_{V,W}(k-1,1) + p_{V,W}(k,1)], \qquad k\in \mathbb{Z},
\end{align*}
and the parameter $q$ is,
%
\be
q := \sum_{k \in \mathbb{Z}}\min\{p(k),p(k+1)\}>0,
\label{qpositive}
\ee
where the positivity of $q$ follows from the fact that
the maximal span is~1.

For the proof we need the following elementary lemma.
It says that, if we 
wait long enough, there will be an (approximately) 
symmetric Bernoulli 
step hidden in $S_n$.

\begin{lemma} \label{Sndecomp}
Under the assumptions of Theorem~\ref{thm:entropycD},
suppose the $X_i$ have zero mean and take values
in $\{a+k\;:\;k\in\IN\}$, for some $a\in\RL$, with maximal span $h=1$.
Then, for each $n \geq 1,$ there is a 
random variable $V^{(n)}$ with values in $\{na+k\;:\;k\in\IN\}$
and a $W^{(n)}\sim{\rm Bern}(q^{(n)})$,
such that,
\begin{equation}\label{sndecompeq}
S_n \equald V^{(n)} + W^{(n)}B,
\end{equation}
where $B\sim{\rm Bern}(1/2)$ 
is independent of $(V^{(n)},W^{(n)})$ and
$q^{(n)}\to 1$ as $n \rightarrow \infty$. Furthermore, 
\begin{equation} 
\label{lemmavarpart}
\mathrm{Var}\big(S_n \big| W^{(n)} = 1 \big) 
= n\sigma^2(1 + o(1)) \quad \text{as } n \rightarrow \infty.
\end{equation}
\end{lemma}

\noindent
{\sc Proof. }
Using the Bernoulli part decomposition $X_i=V_i+W_iB_i$ for
each $X_i$, we can write,
\begin{equation*} 
S_n\equald\sum_{i=1}^n{V_i} + \sum_{i=1}^{N_n}{B_i},
\end{equation*} 
where $N_n = \sum_{i=1}^n{W_i}\sim {\rm Bin}(n,q)$. But also,
\begin{equation*} 
\sum_{i=1}^n{V_i} + \sum_{i=1}^{N_n}{B_i} \equald \sum_{i=1}^n{V_i} 
+ \Biggl(\sum_{i=1}^{N_n - 1}{B_i}\Biggr)\mathbb{I}_{\{N_n \geq 1\}} 
+ \mathbb{I}_{\{N_n \geq 1\}}B,
\end{equation*}
where $B\sim{\rm Bern}(1/2)$ is
independent of everything else. 
This is exactly of the required form~\eqref{sndecompeq},
with $V^{(n)} = \sum_{i=1}^n{V_i} 
+ (\sum_{i=1}^{N_n - 1}{B_i})\mathbb{I}_{\{N_n \geq 1\}}$ 
and $W^{(n)} = \mathbb{I}_{\{N_n \geq 1\}}$,
where $q^{(n)}=1-(1-q)^n\to 1$ as $n\to\infty$
by~(\ref{qpositive}).

For \eqref{lemmavarpart} we only have to consider the case $q < 1,$ 
since otherwise the result holds trivially. 
For the mean we have,
\begin{align*}
0 = \mathbb{E}S_n &= \mathbb{E}\left[\sum_{i=1}^n{V_i}\,\Big|\, W^{(n)} 
= 0\right](1-q)^n
+ \mathbb{E}\bigl(S_n\big| W^{(n)}=1 \bigr)
[1-(1-q)^n].
\end{align*}
On the event $\{W^{(n)} = 0\}=\{W_1=\cdots=W_n=0\}$ the 
$V_i$ are i.i.d., so,
$\mathbb{E}[\sum_{i=1}^n{V_i}| W^{(n)} = 0] = O(n)$, 
and since $[1-(1-q)^n]\to 1$, we must have,
\be
\mathbb{E}\bigl(S_n\big| W^{(n)}=1\bigr) = o(1),
\qquad\mbox{as}\;n\to\infty.
\label{eq:first}
\ee
For the second moment we similarly have,
$$
n\sigma^2 = \mathbb{E}S_n^2 
= \mathbb{E}\left[\left(\sum_{i=1}^n{V_i}\right)^2\,\Big|\, W^{(n)} = 0\right]
(1-q)^n
+ \mathbb{E}\bigl(S_n^2\big| W^{(n)}= 1\bigr)
[1-(1-q)^n],
$$
and since the $V_i$ are i.i.d.\ on $\{W^{(n)}=0\}$,
we have,
$\mathbb{E}[(\sum_{i=1}^n{V_i})^2| W^{(n)}= 0] =O(n^2)$.
Therefore,
\begin{equation}
\mathbb{E}\bigl(S_n^2\big| W^{(n)}=1\bigr) = n\sigma^2(1+o(1)),
\qquad \text{as } n \rightarrow \infty.
\label{eq:second}
\end{equation}
The result follows from~(\ref{eq:first}) and~(\ref{eq:second}).
\qed

We can finally give the proof of the general case of our
main result.

\medskip

\noindent
{\sc Proof of Theorem~\ref{thm:entropycD}. } 
In view of~\eqref{entropyupperbound}, 
we only need to show that, as $n\to\infty$,
\be
H(S_n)\geq\log\frac{\sqrt{n}}{h}+\frac{1}{2}\log(2\pi e\sigma^2)+o(1).
\label{eq:lowerb}
\ee

Without loss of generality, we assume that the $X_i$ have
mean zero and maximal span $h=1$.
Let $\epsilon>0$ be arbitrary and $M$ a large integer to be 
chosen later. For $1\leq i\leq n/M$, 
let $S_i^{(M)}=\sum_{j=(i-1)M+1}^{iM}X_j$, so that
$S_n = \sum_{i=1}^{n/M}S_i^{(M)}.$ 
In the notation of Lemma~\ref{Sndecomp}, for $n \geq M,$
\begin{align*}
H(S_n) 
& = H\Biggl(\sum_{i=1}^{n/M}S_i^{(M)}\Biggr)\\
& = H\Biggl(\sum_{i=1}^{n/M}\Bigl(V_i^{(M)} + W^{(M)}_iB_i \Bigr)\Biggr) \\ 
& \geq  H\Biggl(\sum_{i=1}^{n/M}\Bigl(V_i^{(M)} + W^{(M)}_iB_i \Bigr) 
	\Bigg| W_1^{(M)},\ldots,W_{n/M}^{(M)}\Biggr).
\end{align*}
Let $W^{(M)}$ denote the vector
$(W_1^{(M)},\ldots,W_{n/M}^{(M)})$
and write $A_M$
the collection of vectors $w=(w_1,\ldots,w_{n/M})\in\{0,1\}^{M/n}$ 
with $w_i=1$ for at least $n(q^{(M)}-\epsilon/2)/M$ indices $i$,
where $q^{(M)}=q_i^{(M)}$ is the parameter in
the Bernoulli decomposition of Lemma~\ref{Sndecomp}.
Then we can bound,
\begin{align}
H(S_n)
&\geq \sum_{w \in A_M}
{\mathbb{P}\big(W^{(M)} = w\big)H\Biggl(\sum_{i=1}^{n/M}\Bigl(V_i^{(M)} + W^{(M)}_iB_i \Bigr) \Bigg| W^{(M)} = w\Biggr)}.
\label{averageentropies}
\end{align}
Now observe that, on the event $\{W^{(M)}= w\},$ 
the $n/M$ random variables
$\{V_i^{(M)} + W^{(M)}_iB_i \}$ are independent, though 
not necessarily identically distributed. 
But by \eqref{entropyadd}, we can leave out of the sum inside 
the entropy in \eqref{averageentropies} the summands that correspond 
indices $i$ for which $w_i = 0.$
Thus, writing 
$\bar{W}^{(M)}$ for the vector consisting of $W_i$
with $1\leq i\leq n(q^{(M)}-\epsilon/2)/M$,
and $\mathbf{1}$ for the vector of all 1s,
\begin{align}
H(S_n) &\geq \sum_{w \in A_M} 
{\mathbb{P}\big(W^{(M)} = w\big)
H\Biggl(\sum_{\substack{1 \leq i \leq \frac{n}{M}\;:\;w_i = 1}}
{\Bigl(V_i^{(M)} + W^{(M)}_iB_i \Bigr)} \Bigg| W^{(M)} = w\Biggr)} 
\nonumber\\
&\geq \mathbb{P}\Biggl(\sum_{i=1}^{n/M}{W^{(M)}_i \geq 
\frac{n}{M}\Big(q^{(M)}-\frac{\epsilon}{2}\Big)}\Biggr)
H\Biggl(\sum_{i=1}^{n(q^{(M)}-\epsilon/2)/M}
\Bigl(V_i^{(M)} + B_i \Bigr) \Bigg| \bar{W}^{(M)} = \mathbf{1}\Biggr),
\label{eq:drop}
\end{align}
where the second inequality 
follows form another application of~\eqref{entropyadd},
and the fact that, for different $i$, the distribution of
$V_i^{(M)} + W^{(M)}_iB_i$ only depends on $W^{(M)}_i.$

Since each $W_i^{(M)}\sim{\rm Bern}(q^{(M)})$, the probability
in~(\ref{eq:drop}) converges to~1 exponentially fast.
And since the summands inside the entropy in~(\ref{eq:drop})
are i.i.d.\ with variance $O(1)$, from the upper bound 
in~(\ref{entropyupperbound}) it follows that,
as $n \rightarrow \infty$,
$$H(S_n)
\geq H\Biggl(\sum_{i=1}^{n(q^{(M)}-\epsilon/2)/M}
\Bigl(V_i^{(M)} + B_i \Bigr) \Bigg| 
\bar{W}^{(M)} = \mathbf{1}\Biggr) - o(1).
$$
To complete the proof, we apply Theorem~\ref{bernsmooth} 
to the sequence of i.i.d.\ 
random variables $\{V^{(M)}_i\}$ conditional on 
$\{W^{(M)}_i = 1\}$,
and the independent sequence $\{B_i\}$,
to obtain that, as $n\to\infty$,
$$H(S_n) 
\geq  \frac{1}{2}\log{\Big(\frac{n}{M}
	\Big(q^{(M)}-\frac{\epsilon}{2}\Big)\Big)} 
	+ \frac{1}{2}\log{\Bigl(2\pi e \mathrm{Var}\bigl(V_1^{(M)} 
	+ B_1\big|W_1^{(M)}=1\bigr)\Bigr)} - o(1),
$$
and using the variance bound in Lemma~\ref{Sndecomp},
$$H(S_n) 
\geq 
	\frac{1}{2}\log{\Bigl(\frac{n}{M}\Bigr)} 
	+ \frac{1}{2}\log{\Bigl(2\pi eM\sigma^2\Bigr)} 
	+ \frac{1}{2}\log(1 - \epsilon)
	+ \frac{1}{2}\log{\Big(q^{(M)}-\frac{\epsilon}{2}\Big)} - o(1),
$$
where $M$ is taken large enough for the $o(1)$ term
in Lemma~\ref{Sndecomp} to be smaller than $\epsilon$.
And taking $M$ large enough so that $q^{(M)}>1-\epsilon/2$,
$$
H(S_n)\geq 
	\frac{1}{2}\log{n} + \frac{1}{2}\log{(2\pi e\sigma^2)} 
	+ \log{(1-\epsilon)} - o(1).
$$
Since $\epsilon > 0$ was arbitrary, this gives~(\ref{eq:lowerb})
and completes the proof.

Finally we remark that,
in order to avoid non-essential technicalities,
throughout the proof we have implicitly 
assumed that both $n/M$ and $n(q^{(M)}-\epsilon/2)/M$ are integers.
This does not harm generality as we could have replaced these 
quantities with their integer parts and `` $=$'' with `` $\geq$'' 
where necessary to obtain exactly the same result. 
\qed


\def\cprime{$'$}

\end{document}